\documentclass [11pt,reqno]{amsart}
\usepackage{amsmath, amssymb, amscd}
\usepackage{pb-diagram}
\usepackage{lamsarrow}
\usepackage{pb-lams}
\usepackage{Bdiagx}
\usepackage{amsmath}
\usepackage{amssymb}
\usepackage{yfonts}
\usepackage{color}

\newtheorem{theorem}{Theorem}[section]
\newtheorem{definition}[theorem]{Definition}
\newtheorem{lemma}[theorem]{Lemma}

\newtheorem{corollary}[theorem]{Corollary}

\newtheorem*{remark}{Remark}

\newcommand{\abs}[1]{\lvert#1\rvert}

\renewcommand{\epsilon}{\varepsilon}

\DeclareMathAlphabet{\mathpzc}{OT1}{pzc}{m}{it}
\usepackage{mathrsfs}

\newcommand{\DT}{T^2 \times D^2}

\newcommand{\Z}{\mathbb{Z}}
\newcommand{\C}{\mathbb{C}}

\newcommand{\CP}{\mathbb{CP}}

\newcommand{\R}{\mathbb{R}}
\newcommand{\Slz}{\text{Sl}(3,\Z)}
\newcommand{\Glz}{\text{Gl}(3,\Z)}
\renewcommand{\qed}{$\hfill \square$ \bigskip \\}
\renewcommand{\phi}{\varphi}

\renewcommand{\bar}{\overline}

\address {Laboratoire d'Analyse, de Topologie, et de Probabilit\'es,
Universit\'e de Provence, 13013 Marseille}
\email {zentner@cmi.univ-mrs.fr}

\begin{document}

\title{A note on logarithmic transformations on the Hopf surface}
\author{Raphael Zentner}

\begin{abstract}
In this note we study logarithmic transformations in the sense of differential
topology on two fibers of the Hopf surface. It is known that such
transformations are susceptible to yield exotic smooth structures on
four-manifolds. We will show here that this is not the case for the Hopf
surface,
all integer homology Hopf surfaces we obtain are diffeomorphic to the standard
Hopf surface.
\end{abstract}

\maketitle

\section{Introduction}

The (standard) Hopf surface $S^1 \times S^3$ fibres over the 2-sphere $S^2$ via
the map obtained by composing the   Hopf fibration $S^3 \rightarrow S^2$ with
the projection on the second factor. Any fibre is diffeomorphic to the torus
$T^2$ and there are no singular fibres, because this map is a submersion. It is
a natural problem  to study the effect of logarithmic transformations on two
fibres in this case. Indeed, this operation was successfully used in the case of
the K3 surface to construct exotic K3 surfaces, as well as on other elliptic
fibrations.  These results have been obtained using gauge theoretical methods,
which only apply for manifolds with $b_2^+\geq 1$ \cite{DK} \cite{FM}
\cite{K} \cite{OV}.
Note that all K3-surfaces are diffeomorphic four-manifolds, and there exist
complex
K3-surfaces which are elliptic fibrations. 
In the case of the K3-surface the resulting manifolds depend only on the
multiplicities of the logarithmic transformations, but in our considerations
they depend on some additional parameters as well.

For four-manifolds with the rational  homology of a Hopf surface the existing
gauge
theoretical methods do not apply. On the other hand it is a fundamental and open
problem whether four-manifolds with small second Betti-number 
(especially the four-sphere and
the Hopf surface) do admit exotic structures. The four-manifold with
smallest second Betti number admitting exotic smooth structures which is
known at present is $\mathbb{CP}^2 \# 5 \bar{\mathbb{CP}}^2$ \cite{PSS}. In the 
complex geometric
framework, exotic Hopf surfaces do not exist, for by a result of Kodaira
\cite{Ko} every complex surface which is homeomorphic to $S^1\times S^3$ is a
primary Hopf surface, so it is diffeomorphic to $S^1\times S^3$.  Complex
surfaces which are rational homology Hopf surfaces have been classified in
\cite{EO} using logarithmic transformations. Further results about elliptic
surfaces in the class of complex surfaces can be found in \cite{FM}. Our
situation here, however,
is purely topological in nature, and the logarithmic transformations considered
are more general than the complex-geometric ones. In particular, logarithmic
transformations with multiplicity zero do not arise in  the complex geometric
setting, and may even result in manifolds not admitting any complex structure at
all \cite{G}.

We will first calculate the fundamental group of the manifold obtained by two
logarithmic transformations. As it turns out in many cases, including
multiplicity 0, the resulting manifold will have the same fundamental group as
the Hopf surface. Since the Euler characteristic is invariant under logarithmic
transformations, we will   obtain a manifold having the same (integer) homology
as the Hopf surface. We will then describe a procedure to construct all these
manifolds by gluing two copies of $T^2 \times D^2$ via a diffeomorphism between
their boundaries. Using diffeomorphisms of $T^2 \times S^1$ which extend over
$T^2 \times D^2$, we will be able to show that manifolds given by different
gluing diffeomorphisms may still be diffeomorphic. Using this observation, we
will  find  a certain standard form for every homology Hopf surface obtained by
this gluing method. The possible standard forms are determined by elements in
$\text{Sl}(2,\Z)$. Finally, using a handlebody-theoretical argument \cite{LP},
we  prove that this parameter does not affect the diffeomorphism type. \\
%
%
%
\section*{Acknowledgements}

I am grateful to Peter Kronheimer for helpful conversations on this and related
topics. I am also indebted to my advisor Andrei Teleman for proof-reading the
paper and related suggestions, as well as for the encouragement to write this
paper. Furthermore I am thankful to the referee for some useful
comments. Finally, I wish to thank Amy Ellingson for proof-reading the
English.
\section {Logarithmic transformations applied to Hopf surfaces and resulting
fundamental group}
\begin{definition} Let $\pi : X \rightarrow \Sigma$ be an elliptic fibration. We
say that a four-manifold $X'$ is obtained from $X$ by logarithmic transformation
on
a regular fibre $F$ of $\pi$ if   $X'$ is obtained from $X$ through the
following construction: We cut out a regular neighbourhood $\nu F$ of $F$ and we
glue in a  $\DT$ via an arbitrary orientation-reversing  diffeomorphism
$\varphi: T^2 \times S^1 \rightarrow \partial \nu F$. The absolute value of the
degree of $\pi|_{\partial \nu F} \circ \varphi|_{pt \times S^1}$ is called the
multiplicity of the logarithmic transformation \cite{G}.
\end{definition}

The diffeomorphism $\varphi$ is determined, up to isotopy, by its induced
isomorphism of fundamental groups, which is itself, after the choice of some
bases,
is determined by a matrix in $\Glz$. Alternatively, we fix one such
diffeomorphism, which can be used to identify $\partial \nu F$ with $T^2 \times
S^1$. Any other diffeomorphism is determined by a self-diffeomorphism of $T^2
\times S^1$,
and these diffeomorphisms are given, up to isotopy, by elements in $\Slz$.

We will first give a gluing description of the Hopf-surface $X = S^1 \times S^3$
which will turn out useful. For this we shall first describe
$S^3$ as two solid tori $S^1 \times D^2$ glued together. The two closed discs
$D^2$ will turn out to be the northern and southern hemisphere,
respectively, under
the Hopf fibration $S^3 \to S^2$. Indeed, $S^3$ can be seen as the following
set:
\[
  S^3 = \left\{ \left. (z,w) \in \C^2  \right| \ \ \abs{z}^2+\abs{w}^2 = 2
\right\}
\]
The Hopf fibration is then given by the map 
\[
  S^3 \to \CP^1 \quad \text{given by} \quad (z,w) \mapsto [z:w] \ ,
\]
and $\CP^1$ is diffeomorphic to $S^2$. Define $S^3_+$ to be the set of elements
$(z,w)$ such that $0 \leq \abs{w}^2 \leq 1$, and $S^3_-$ to be the set of
elements $(z,w)$ with $0 \leq \abs{z}^2 \leq 1$. Then there are diffeomorphisms 
\begin{equation*}
\begin{split}
S^3_+ & \xrightarrow{f_+} S^1 \times D^2, \quad \text{given by} \quad 
  (z,w) \mapsto \left( \frac{z}{\abs{z}} , \frac{w}{z} \right) , \quad
\text{and} \\
S^3_- & \xrightarrow{f_-} S^1 \times D^2, \quad \text{given by} \quad 
  (z,w) \mapsto \left( \frac{w}{\abs{w}} , \frac{z}{w} \right)  .
\end{split}
\end{equation*}
When we restrict $f_+ \circ f^{-1}_-$ to the boundary, then the map $\partial
(S^1 \times D^2) \to \partial (S^1 \times D^2)$ is given by the formula
\begin{equation*}
  f_+ \circ f_-^{-1} (u,\xi) = (u \xi , \overline{\xi}) \ .
\end{equation*}
We extend this latter map to the trivial $S^1$ factor by the identity, so that
we get a map $\zeta: T^2 \times \partial D^2 \to T^2 \times \partial D^2$, and
the description of the Hopf surface as a gluing
\begin{equation}\label{Hopf}
  X = (T^2 \times D^2) \cup_{\zeta} (T^2 \times D^2) \ .
\end{equation}
Now let us consider the manifold $X'$ obtained from the Hopf surface when
performing logarithmic transformations on two fibres, say on the fibre $F_+$
over the north pole $x_+:=[1:0]$ and the fibre $F_-$ over the south pole
$x_-=[0:1]$, associated with diffeomorphisms $\varphi_\pm$. There are natural
identifications of $\partial (X - \nu F_\pm)$ with the "inner" boundary of $T^2
\times (D^2 - \mathring{D}^2_{1/2})$ according to the decomposition
(\ref{Hopf}). Therefore the orientation-reversing diffeomorphisms $\varphi_\pm$
can be seen as an orientation-preserving diffeomorphism of $T^2 \times S^1$,
because the above "inner" boundary is with opposite orientation to the "outer"
boundary. Let us denote by $X_\pm$ the two manifolds $(T^2 \times (D^2 -
\mathring{D}^2_{1/2})) \cup_{\varphi_\pm} (T^2 \times D^2)$. Gluing
two manifolds along their boundaries is actually a suitable
identification of  collar neighbourhoods of the boundaries of the two manifolds.
In our case, this description is given as
\[
 X_{\pm}= \left(T^2 \times \left(D^2 - D^2_{1/3}\right)\right) \cup_{\Phi_\pm}
\left(T^2 \times \mathring{D}^2_{2/3}\right) \ ,
\]
where $
\Phi_\pm: \left(\frac{1}{3},\frac{2}{3}\right) \times T^2 \times S^1 \to
\left(\frac{1}{3},\frac{2}{3}\right) \times T^2 \times S^1$ is given by
$\Phi_\pm\left(r,u,v,\xi \right):=
\left(\frac{2}{3}-r,\varphi_\pm\left(u,v,\xi\right)\right) \ .
$ 
Let us now fix some paths inside $D^2 \times T^2$, where the disc is thought of
as a subset of $\C$, centred at the origin. Fix some base-point $(u_0,v_0,\xi_0)
\in T^2 \times D^2$, where $\abs{\xi_0}=\frac{1}{2}$, so that the base point is
in the ``gluing area". Let us define three paths $\alpha_\pm, \beta_\pm,
\gamma_\pm$ by the formulae 
$\alpha_\pm(t)=(u_0,v_0e^{it},\xi_0), \ \beta_\pm(t)=(u_0 e^{it},v_0, \xi_0), \
\text{and} \ \gamma_\pm(t) = (u_0,v_0, \xi_0 e^{it})$. The path $\gamma_\pm$ is
then a meridian to the fibre $T^2 \times \{ 0 \}$ over $x_\pm$ - its
projection onto the fibre is trivial - whereas $\alpha_\pm$ and $\beta_\pm$
induce a basis of the fundamental group of the fibre. Note that by the same
formulae we can define paths $(\alpha_\pm', \beta_\pm', \gamma_\pm')$ inside the
pieces $T^2 \times D^2$ to be glued in with $\varphi_\pm$. Then
$(\alpha_\pm,\beta_\pm,\gamma_\pm)$ induces a basis of $\pi_1(X - \nu F_\pm)$
and
$ (\alpha_\pm', \beta_\pm', \gamma_\pm')$ induces a basis of $\partial(T^2
\times
D^2)$. The diffeomorphisms $\varphi_\pm$ are then determined by their maps of
fundamental groups 
\[
 \varphi_{*}^+ = \begin{pmatrix} * \ & * \ & a \cr *\ & *\ & b \cr *\ & *\ & p
\end{pmatrix} \quad 
 \varphi_{*}^- = \begin{pmatrix} * \ & * \ & c \cr *\ & *\ & d \cr *\ & *\ & q
\end{pmatrix} \ ,
\]
which are elements in $\Slz$. The entries marked as $*$ will not be relevant to
the fundamental group, as we shall see. We call $(a,b) \in \Z^2$ the
direction of the logarithmic transformation $\varphi_+$, and $\abs{p}$ is its
multiplicity.

In order to compute the fundamental group of $X'$ we shall first compute the
fundamental groups of $X_\pm$ and then glue them together via $\zeta$. $X_+$ is
given as the union of two open sets, namely the sets $X_1 = T^2 \times (D^2
-D^2_{1/3})$ and $X_2 = T^2 \times \mathring{D}^2_{2/3}$, with intersection $X_0
= T^2 \times (\mathring{D}^2_{2/3} - D^2_{1/3})$. The manifold $X_0$ injects
into $X_1$
via the natural inclusion $i$, and into $X_2$ via $\phi_+$. The fundamental
group of each piece is
\begin{equation*}
\begin{split}
 \pi_1(X_0) =&   \langle \alpha_0, \beta_0, \gamma_0 \, |  \, [ \ , \ ]=1
\rangle \ , \\
 \pi_1(X_1) = & \langle \alpha, \beta, \gamma \, |  \, [ \ , \ ]=1 \rangle \ ,
\quad \text{and} \\
 \pi_1(X_2) = & \langle \alpha', \beta' \, |  \, [ \ , \ ]=1 \rangle \ .	
\end{split}
\end{equation*}
By $[ \ , \ ]$ we simply mean that all commutator relations are satisfied. The
Seifert-van Kampen theorem states that $\pi_1(X_+)$ has as generators
both the generators of $\pi_1(X_1)$ and of $\pi_1(X_2)$, as relations all
those of $\pi_1(X_1)$
and $\pi_1(X_2)$, and the additional relations 
\begin{equation*}
\begin{split}
i(\alpha_0) = & \phi(\alpha_0) \ \Leftrightarrow \ \alpha' = \phi(\alpha_0), \\
i(\beta_0)  = & \phi(\beta_0) \ \Leftrightarrow \ \beta' = \phi(\beta_0), \quad
\text{and} \\
i(\gamma_0) = & \phi(\gamma_0) \ \Leftrightarrow \ 1 = \phi(\gamma_0) \ .
\end{split}
\end{equation*}
The first two relations imply that we can drop the generators $\alpha'$ and
$\beta'$ as well as these two relations. Therefore the fundamental group is
\begin{equation*}
\pi_1 (X_+) = \langle \, \alpha_+, \beta_+, \gamma_+ \,|\,  [\ ,\ ]=1,\,
\alpha_+^a \beta_+^b \gamma_+^p=1 \ \rangle \ .
\end{equation*}
Correspondingly, we get 
\begin{equation*}
\pi_1 (X_-) = \langle \, \alpha_-, \beta_-, \gamma_- \,|\,  [\ ,\ ]=1,\,
\alpha_-^c \beta_-^d \gamma_-^q=1 \ \rangle \ .
\end{equation*}
In order to compute the fundamental group of $X'=X_+ \cup_\zeta X_-$ we proceed
in the same way. $T^2$ times a ``middle annulus" injects into $X_-$ via the
natural inclusion, whereas it injects into $X_+$ via $\zeta$. As we have
$\zeta_*(\alpha_0) = \alpha_+, \zeta_*(\beta_0)= \beta_+$ and $\zeta_*(\gamma_0)
= \alpha_+ \gamma_+^{-1}$ we get a final formula:
\begin{equation*}
\pi_1 (X') = \langle \, \alpha, \beta, \gamma \,|\,  [\ ,\ ]=1, 
\alpha^a \beta^b (\alpha \gamma^{-1})^p = 1, \, \alpha^c \beta^d  \gamma^q = 1 \
\rangle \ .
\end{equation*}
By the classification of finitely generated Abelian groups, we find that we have
an isomorphism $\pi_1(X') \cong \Z \oplus \Z/{\mu}\Z $,
where $\mu$ is the highest common divisor of all the 2-minors of a presentation
matrix for this group. It is easy to see that there are various choices possible
in which this number equals $1$, including cases where one or both of the
multiplicities may be zero. 

\begin{remark} If we perform the two logarithmic transformations such that they
are trivial on the $S^1$-factor, then the construction is $S^1$ times
Dehn-surgery on the Hopf-link in $S^3$. The resulting four-manifold is then
$S^1$
times a lens space; this can be seen using the surgery description of lens
spaces \cite{GS}. However, we shall point out that if the logarithmic
transformations are of type $(1 \ , 0 \ , p)$ and $(1 \ , 0 \ , q)$,
then the Dehn-surgery description of the Lens space we get is not the
surgery description on the Hopf link with surgery coefficients $p$ and $q$,
with respect to the blackboard framing. \end{remark}

\section {Formulation in terms of gluing two copies of $T^2 \times D^2$}
We will denote by $ X_\varphi := (\DT) \cup_\varphi( \DT)$ the four-manifold
obtained by gluing $\DT$ to $\DT$ via the orientation-reversing diffeomorphism
$\psi$ between their boundaries. Let us denote by $A^2$ an annulus. There are
canonical identifications of the boundary-components of $T^2 \times A^2$ with
$T^2 \times S^1$, as before.

Two isotopic diffeomorphisms induce diffeomorphic manifolds, so we are only
interested in isotopy classes of diffeomorphisms here. Furthermore we may
restrict our attention to orientation-reversing diffeomorphisms. We shall also
identify the boundaries of the two copies of $\DT$ with the 3-torus $T^2 \times
S^1$, once orientation-preserving, once orientation-reversing, and this once and
for all. The diffeomorphism $\varphi$ is then given by an orientation-{\it
preserving} diffeomorphism of $T^3$. Finally, in the case of the 3-torus, such a
diffeomorphism up to isotopy is determined by its associated automorphism of the
fundamental group, and therefore by a matrix in $\Slz$.
We will show here that all of the manifolds considered so far can be obtained by
gluing just two copies of $T^2 \times D^2$ along their boundaries:
\begin{lemma} We have the following diffeomorphism:
\[
X_{\psi\ \circ \ \varphi} \cong (\DT )\cup_\psi (T^2 \times A^2) \cup_\varphi
(\DT) .
\]
\end{lemma}
{\em Proof:} As any diffeomorphism of one boundary-component of $T^2 \times A^2$
extends over the whole of $T^2 \times A^2$ the result follows easily. \qed

Our next objective is to calculate the fundamental group of $X_\varphi$.
Let us use the bases $(\alpha_\pm, \beta_\pm, \pm \gamma_\pm^{\pm1})$ from above
(up to "orientation") and suppose that the map $\varphi_*$, which is now given
by an element of $\Slz$, looks as follows:
\begin{equation}\label{gluingmatrix}
\varphi_{*} = \begin{pmatrix} a \ & c \ & g \cr b\ & d\ & h \cr e\ & f\ & k 
\end{pmatrix} \quad 
\end{equation}
According to the theorem of Seifert-van Kampen a presentation of the fundamental
group of
$X_\varphi$ is given by
\begin{equation*}
\pi_1 (X_\varphi) = \langle \alpha, \beta \ | \ [\alpha,\beta]=1, \,
(\alpha^{g'} \beta^{h'})^{(g,h)}=1 \ \rangle \ .
\end{equation*}
Here $(g,h)$ denotes the greatest common divisor of $g$ and $h$, and $g'$
and$ h'$ are
such that $g=(g,h)\ g'$, $h=(g,h)\ h'$. We define $(0,0):=0$.
The fundamental group is therefore isomorphic to $\Z \oplus \Z/{(g,h)\Z} $:
\begin{equation*}
\pi_1 (X_\varphi) = \Z \oplus \Z/{(g,h)\Z} . 
\end{equation*}
In particular, $X_\varphi$ is a homology Hopf surface if and only if $(g,h)=1$,
noting that any $X_\varphi$ has Euler-characteristic zero. 

If now we  perform the logarithmic transformations associated with $\varphi_\pm$
on the two fibres $F_\pm$ of the Hopf surface, then the  resulting manifold will
be given by the following gluing construction
\begin{equation*}
(\DT) \, \cup_{\varphi_+^{-1}} (T^2 \times A^2) \, \cup_\zeta \, (T^2 \times
A^2) \, \cup_{\varphi_-} (\DT) 
\end{equation*}
which is diffeomorphic, by the above lemma, to 
\begin{equation*}
X_{\varphi_+^{-1} \circ \ \zeta \ \circ \ \varphi_-} \ .\quad
\end{equation*}
Whether this manifold is a homology Hopf surface can now be deciphered from the
automorphism $\left(\varphi_+^{-1} \circ \ \zeta \ \circ \ \varphi_-\right)_*$
of the fundamental group. However, calculating the entity
$(g,h)$, which a posteriori depends only on the numbers $a,
b, p$ and $c, d, q$, using this matrix product, is a rather difficult
problem.

\begin{theorem}
Suppose the manifold $X_\varphi$, constructed as above, is a homology Hopf
surface. Then $X_\varphi$ is diffeomorphic to the Hopf surface $X_\zeta$.
\end{theorem} 

\begin{corollary}
If logarithmic transformations on two fibres yield a homology Hopf surface then
this four-manifold is diffeomorphic to the standard Hopf surface $S^1 \times
S^3$.
\end{corollary}
 
{\it Proof of the Theorem.} Observe first that the  two manifolds
\begin{equation*}
X_{\varphi} \ ,\   X_{\psi_t^{-1}  \circ \ \varphi \ \circ \ \psi_b}
\end{equation*}
are diffeomorphic as  soon as the diffeomorphisms $\psi_t$ and $\psi_b$ of $T^2
\times S^1$   extend over $T^2 \times D^2$ as diffeomorphisms. A diffeomorphism
$\psi$ extends if and only if the associated matrix has the form 

\begin{equation}\label{ext_matrix}
\psi_* = \left(\begin{array}{ccc} r \ & t\  & 0\  \\ s\  & u\  & 0\  \\ v\  & w\
 & 1\  \end{array} \right) .
\end{equation}
Indeed, it is easy to construct explicitly extensions of these diffeomorphisms;
on the other hand, if $\psi$ extends to a diffeomorphism$\Psi$, then the first
two entries
in the third column of the corresponding matrix must be zero. This can be seen
using  the commutative diagram
\begin{equation*}
\begin{split}
\begin{diagram}\dgmag{800} \dgsquash[3/4]
\node{H_1(T^2 \times S^1)} \arrow[2]{e,t}{\psi_*} \arrow[1]{s} \node[2]{H_1(T^2
\times S^1)} \arrow[1]{s} \\
\node{H_1(T^2 \times D^2)} \arrow[2]{e,t}{\Psi_*} \node[2]{H_1(T^2 \times D^2)
.}
\end{diagram}
\end{split}
\end{equation*}
This observation can be used to perform certain line operations on $\varphi_*$
by
left-multiplication with matrices induced by extending diffeomorphisms, as well
as to perform certain column operations by right-multiplication with these
matrices, without changing the diffeomorphism type.

Suppose now that $X_\varphi$ is a homology Hopf surface with associated matrix
$\varphi_*$ as in (\ref{gluingmatrix}) above.
In particular, the greatest common divisor of $g$ and $h$ is one: $(g,h)=1$. By
left-multiplying with a matrix $u \in\text{Sl}(2,\Z) \subseteq \text{Sl}(3,\Z)$,
where the inclusion is into the upper left part in the $3 \times 3$ matrix, we
may
assume that $g=1, h=0$ in (\ref{gluingmatrix}). Such a matrix $u$ is of type
(\ref{ext_matrix}). Now there is a matrix $L$ of type (\ref{ext_matrix}) such
that left-multiplication of the new matrix $\varphi_*$ by $L$ adds $-(k-1)$
times the first line of $\varphi_*$ to its last line. Therefore we may suppose
that $k=1$. Now there is a matrix $R$ of the type (\ref{ext_matrix}) such that
right-multiplication of the newest $\varphi_*$ by $R$ will add appropriate
multiples of the third column of $\varphi_*$ to its first and second, so that we
may assume $e=f=0$ because $k=1$. $\varphi_*$ in (\ref{gluingmatrix}) may
therefore be supposed to have the form
\begin{equation}\label{end_form}
 \varphi_* = \begin{pmatrix} a\  & c\  & 1\  \cr b\  & d\  & 0\  \cr 0\  & 0\  &
1\  \end{pmatrix} \ .
\end{equation}
A corresponding diffeomorphism is given by $\varphi(u,v,z)=(u^a v^c z, u^b v^d,
z)$.
We can't simplify much further in order to obtain the matrix $\zeta_*$,
where $\zeta$ is inducing the standard Hopf surface as above. However, the
attachment of $T^2 \times D^2$ to the upper $T^2 \times D^2$, which we shall
denote by $X_+$, may be done by attaching first a 2-handle, then two 3-handles,
and eventually a 4-handle. To be more precise, decompose the torus $T^2$ in the
obvious way into a $0$-handle $\Sigma_0$, two $1$-handles $\Sigma_{11}$ and
$\Sigma_{12}$, and a 2-handle $\Sigma_2$. Then the attachment, via $\varphi$, of
$\Sigma_0 \times D^2$ to $X_+$ is done along $\Sigma_0 \times \partial D^2$,
thus
we attach a 2-handle and get $X^{(2)}:=X_+ \cup (\Sigma_0 \times \partial D^2)$.
It is now easily verified that $\Sigma_{11} \times D^2$ and $\Sigma_{12} \times
D^2$ are attached to $X^{(2)}$ along a thickened 2-sphere $S^2 \times
D^1$, corresponding to 3-handle attachment. Finally $\Sigma_2 \times
D^2$ is glued to the resulting manifold along a 3-sphere, a 4-handle
attachment.
The union of the two 3- and the one 4-handle is diffeomorphic to a boundary sum
$S^1 \times D^3 \ \natural \ S^1 \times D^3$, which is the gluing of two pieces
of $S^1 \times D^3$ via a diffeomorphism between two discs in their boundaries.
The boundary of this manifold is $S^1 \times S^2 \ \# \ S^1 \times S^2$.
 It is known \cite{LP} that any diffeomorphism of $S^1 \times S^2 \ \# \ S^1
\times S^2$ extends over the whole boundary sum. Therefore only the
2-handle-attachment is relevant for determining the diffeomorphism type of the
closed four-manifold.
  
On the other hand, the attaching of $\Sigma_0 \times \partial D$ is determined,
up to isotopy, by the attachment of the attaching sphere $\{0\} \times S^1$ as
well as the isomorphism of normal bundles $\nu_{\Sigma_0\times S^1}(\{0\} \times
S^1) \to \nu_{T^3}(\varphi(\{0\}\times S^1))$ induced by the derivative
$d\varphi$. We shall denote by $L_\varphi$ this bundle isomorphism. After
identification of $\Sigma_0$ with a ball centred in the origin in $\R^2$ we get
a canonical isomorphism $\nu_{\Sigma_0\times S^1}(\{0\} \times S^1) \cong S^1
\times \R^2$. By a framing $f$ of $\varphi(\{0\} \times S^1)$ we understand a
fixed isomorphism of the normal bundle $\nu_{\Sigma_0\times S^1}(\{0\} \times
S^1)$ with $S^1 \times \R^2$. We say that a framing $f$ is isotopic to the
framing $f'$ if they are homotopic through bundle isomorphisms. By replacing
$L_\varphi$ with $f^{-1}$ we see that the 2-handle attachment is determined by
$(\varphi(\{0\}\times S^1),f)$, the embedding with a given framing of the
attaching sphere. Thus, framings and the isomorphisms $L_\varphi$ are
equivalent
notions.
Up to isotopy, the attachment depends only on the framing up to isotopy. If we
fix one framing, we see that all possible isomorphisms of normal bundles are
given by bundle automorphisms of $S^1 \times \R^2$. 


For the above choice of $\varphi$ the attachment of the attaching sphere does
not depend on 
the specific entries in $\varphi_*$. We identify the normal bundle of
$\varphi(\{0\}\times S^1)$ with orthogonal complement to its tangent bundle
within $T(T^3)$, and get an identification with $S^1 \times \R^2$ by specifying
two constant orthonormal sections of that bundle,
$ e_1 = (1 \, , \, 0 \, , \, -1)$ and $e_2= ( 0 \, , \, 1 \, , \, 0)$.
The isomorphism $L_\varphi$ is then given by the {\em constant} matrix
\[
L_\varphi = \begin{pmatrix} a & c \\ b & d \end{pmatrix} \ .
\]
Because this matrix is in $\text{Sl}(2,\Z)$ we see that there is an isotopy of
bundle automorphisms taking one automorphism into the other. In other words, the
corresponding
framings are isotopic. \qed

\end{document}